\documentclass[12pt,a4paper]{article}
\usepackage{euscript,amsfonts,amssymb,amsmath,amscd}

\input{epsf}
\newcommand{\op}{\operatorname}
\newcommand{\m}{\mathbb}
\newcommand{\e}{\EuScript}
\begin{document}

\title{Complex cobordisms and singular manifolds arising from Chern classes.}

\author{A.Kustarev.}

\maketitle

\begin{abstract}
This paper deals with the question of J.Morava on existence of canonical
complex cobordism class of singular submanifold. We present 
several solutions of this question for $X_r(\xi)$ -- the set of points where
$\dim\xi-r+1$ generic sections of a complex vector bundle $\xi$ are linearly
dependent. The corresponding complex cobordism classes $Q_r(\xi)$ and $P_r(\xi)$
tend to have many nice properties, such as deformed sum formula,
but they don't coincide with Chern classes $c_r^U(\xi)$. 
They also have relation to the theory of $IH$-small
resolutions.

\end{abstract}

\section{Introduction.}

A well-known question of Steenrod (\cite{steenrod}) asks if a given homology
class in some 
cellular complex may be realized as an image of fundamental class of
some manifold. As shown by R.Thom (\cite{thom}), the Steenrod problem always has
the positive solution in $\m Z_2$-homology. But if the ring of coefficients 
is $\m Z$, there are counterexamples to the problem, though one can still show
that an arbitrary class taken with large multiplicity may be realized as the
image of fundamental class of a manifold. 

Let $Y\subset M$ be a complex semialgebraic subset, that is, locally determined
by algebraic equations, in a manifold $M$. (In this paper $M$ is compact,
even-dimensional, oriented manifold of real dimension $2m$ 
with no boundary). Then the set of 
singular points of $Y$ has real codimension at least two in $Y$, so the homology
(and dual cohomology) class $[Y]\in H_*(M,\m Z)$ is well-defined.\\  

{\bf Question (J. Morava).} Does there exist a {\it canonical} complex cobordism
class $[Y]\in U^*(M)$ related to $Y$?\\

The complex cobordism ring $U^*(M)$ is 
generated by "singular bordisms" -- maps of smooth manifolds $X\to M$ with
complex structure in stable normal bundle (for more details see \cite{KF} or
\cite{stong}). Clearly one can find at least one resolution of $Y$ by applying 
Hironaka's theorem on singular varieties, but usually it gives
very uneffective solution of the problem. One has to look for efficient canonic
resolutions depending on definition of $Y$.

In the present paper we deal with the following situation. Let $\xi$ be an
$n$-dimensional complex vector bundle over $M$. One can consider the set
$X_r(\xi)\subset M$ where generic $n-r+1$ sections of $\xi$ are
linearly dependent. Then, according to well-known definition 
going to Pontryagin (\cite{pontryagin}), we have 
$[X_r(\xi)]=c_r(\xi)$ in integral
cohomology. 
In this paper we give some explicit solutions of Morava question for
$Y=X_r(\xi)$ and investigate the corresponding cobordism classes. 

We begin with description of several 
geometric constructions of cohomological 
Chern classes, based on well-known constructions in singularity theory and Thom
polynomials. Here "geometric" means "representing dual homology class as an
image of fundamental class of some canonically defined manifold". It is
known that top Chern class (Euler class) can always be represented by even an
embedded submanifold (as follows from Thom transversality theorem). Using only
notions of Euler class and pushforward (Gysin) homomorphism we show that
Steenrod problem for $c_r(\xi)$ always has the positive solution. (So there are
examples of integer cohomology classes which can't be Chern classes of any
complex vector bundle).

Our first result (Prop. 3.3) says that Chern classes in complex cobordisms
$c_r^U(\xi)$ 
don't solve Morava question for $X_r(\xi)$. (The classes $c_r^U(\xi)\in U^{2r}$
are uniquely determined by the same four axioms, see \cite{KF}). So one should
look for some other cobordism characteristic classes, which are represented by
maps $X\to X_r(\xi)$ resolving the singularities of $X_r(\xi)$. These classes
are named by $Q_r(\xi)$ and $P_r(\xi)$, we investigate their basic properties
(Th. 4.2 and 7.2) -- they satisfy many Chern classes axioms, except Whitney sum
formula and triviality in non-positive dimensions. 

The deformation of sum formula for classes $Q_r(\xi)$ may be found explictly
(Th. 6.3), thus giving the possibility of their purely axiomatic definition
(Th. 6.4).

The classes $Q_r(\xi)$ and $P_r(\xi)$ are different (Th. 7.3), though in case of
small dimensions of $M$ (up to $\dim_{\m R} = 8$) they coincide. If we switch to
algebraic situation ($M$ is a nonsingular variety, $\xi$ is a locally free
sheaf), then $X_r(\xi)$ is a singular variety. Our resolutions turn to be
$IH$-small resolutions of $X_r(\xi)$ (see \cite{totaro}) and, as shown in that
paper, their complex elliptic genera are equal. That gives restrictions on
$P_r(\xi)-Q_r(\xi)$. 

The author is grateful to his advisor, Prof. V.M.Buchstaber, for attention to
his work and to Prof. M.E.Kazarian for valuable discussions.

\section{Geometric constructions of cohomological Chern classes.}

In this section we present some new definitions of cohomological
Chern classes. These
definitions use only notions of Euler class and pushforward homomorphism, so
the Poincare duals of $c_r(\xi)$ are represented as the images of fundamental
classes of nonsingular oriented manifolds, thus giving the solution of Steenrod
problem for $c_r(\xi)$. 

Let $\xi$ be an $n$-dimensional complex vector bundle over a manifold $M$. A well-known definition of Chern classes 
(going to Pontryagin, \cite{pontryagin}) deals with set $X_r(\xi)$ -- a set of
points on $M$ where $(n-r+1)$ generic sections of $\xi$ are not of maximal rank.
The set $X_r(\xi)$ is a complex semialgebraic set, that is, locally determined
by algebraic equations, so its homology and dual cohomology classes are
well-defined since all the singularities have real codimension at least two. 
The cohomology class $[X_r(\xi)]$ is, by definition, the Chern class 
$c_r(\xi)$.

In the further text we'll consider a smooth generic section $s$ of a bundle
$\op{Hom}(\m C^{n-r+1},\xi)$ instead of $(n-r+1)$ sections of $\xi$. The Steenrod problem always has positive solution for classes 
$c_r(\xi)$ -- we are going to present some explicit constructions starting from
the $X_r(\xi)$-definition of $c_r(\xi)$. They may also be regarded as 
the alternate definitions of Chern classes.\\

{\bf Proposition 2.1.} {\it Let $\gamma$ be a tautological bundle over 
projectivization
$M\times \m CP^{n-r}$ of a trivial $(n-r+1)$-dimensional bundle $\m C^{n-r+1}$ 
and $p:M\times \m CP^{n-r}\to M$ -- a projection map. 
Then $r$-th Chern class $c_r(\xi)$ 
of bundle $\xi$ equals to $p_! e(\gamma^*\otimes p^*\xi)$.}\\

{\it Proof 1.} The $(n-r+1)$ generic sections of $\xi$ over $M$ define the 
section $s_Q$ 
of bundle $\op{Hom}(\gamma,p^*\xi)$ over $M\times\m CP^{n-r}$ by taking
the composition: $\gamma\to \m C^{n-r+1}\to p^*\xi$. The set $X_r^Q(\xi)$ 
of points on $M\times \m CP^{n-r}$ where $s_Q$ vanishes is a nonsingular
submanifold by Thom transversality theorem and its cohomology class
coincides with $e(\op{Hom}(\gamma,p^*\xi))=c_n(\op{Hom}(\gamma,p^*\xi))$. 

Note that $p(X_r^Q(\xi))=X_r(\xi)$ and if rank of $s$ drops by $l$ in point
$x\in M$, then $p^{-1}(x)\cap X_r^Q(\xi) = \m CP^{l-1}$. 
The generic point $x\in X_r(\xi)$ has
 $l=1$ and $p^{-1}(x)\cap X_r^Q(\xi)$ 
is just a single point. So $X_r^Q(\xi)$ is a nonsingular
resolution of $X_r(\xi)$ and $p_*([X_r^Q(\xi)])=[X_r(\xi)]$. $\Box$

{\it Proof 2.} We can calculate $p_! e(\gamma^*\otimes p^*\xi)$ explicitly using
splitting principle. If $t=c_1(\gamma^*)$ and $t_1\dots t_n$ are 
Chern roots of $\xi$, so
$c(\xi)=(1+t_1)\dots(1+t_n)$, then $e(\gamma^*\otimes 
p^*\xi)=(t+t_1)\dots(t+t_n)$. Pushforward homomorphism $p_!:H^*(M\times \m
CP^{n-r})\to H^*(M)$ acts by formula $p_!(t^{n-r})=1$, $p_!(t^i)=0, i\ne n-r$.
It follows that $p_! e(\gamma^*\otimes p^*\xi)$ is $r$-th elementary symmetric
polynomial in $t_1\dots t_n$, which is $c_r(\xi)$. $\Box$

Now we define the "dual" resolution $X_r^P(\xi)$ of $X_r(\xi)$. It corresponds
to section $s^*:\xi^*\to \m C^{n-r+1}$ dual to $s:\m C^{n-r+1}\to\xi$.\\

{\bf Proposition 2.2.} {\it Let $\gamma$ be a tautological $r$-dimensional 
bundle over the grassmanization $Gr_r(\xi^*)$ of $\xi^*$ 
and $p:Gr_r(\xi^*)\to M$ -- 
a projection map. Then the Chern class $c_r(\xi)$ equals to $p_!
e((\gamma^*)^{n-r+1})$.}\\

The proof is analogous to the proof 1 of Prop. 2.1. The section $s^*$ of
$\op{Hom}(\xi,\m C^{n-r+1})$ determines the section $s^*_P$ of
$\op{Hom}(\gamma,\m C^{n-r+1})$. If the rank of $s$ in the point $x\in M$ 
drops by 
$l$, then the rank of $s^*$ drops by $r+l-1$ in $x$. If $x$ is a generic point
of $X_r(\xi)$, the rank in $x$ drops by $r$ and the preimage $p^{-1}(x)\cap
X_r^P(\xi)$ is a
single point. The set $X_r^P(\xi)$ of
points on $Gr_r(\xi^*)$ where $s^*_P$ vanishes is a nonsingular resolution of
$X_r(\xi)$, so the image of a fundamental class $p_*([X_r^P(\xi)])$ equals to
$[X_r(\xi)]$. $\Box$

If $r>1$, then $X_r^Q(\xi)$ and $X^P_r(\xi)$ are different resolutions of
$X_r(\xi)$. The preimage of point $x\in X_r(\xi)$ where $\dim\ker s = l$ is $\m
CP^{l-1}$ in $X_r^Q(\xi)$ and $Gr_r(\m C^{r+l-1})$ in $X_r^P(\xi)$. If $r=1$,
then $\m CP^{l-1} = Gr_r(\m C^{r+l-1})$ but the corresponding resolutions are
different -- their complex cobordism classes differ as we show later. 

One can also define the Chern classes in cohomology using the transfer
homomorphism:\\

{\bf Proposition 2.3 (\cite{viniti})} {\it Let $\gamma$ be a $r$-dimensional 
tautological bundle over grassmanization
$Gr_r(\xi)$ and $t:H^*(Gr_r(\xi))\to H^*(M)$ a transfer homomorphism. 
Then $t(e(\gamma)) = c_r(\xi)$.}\\

Cohomological transfer homomophism $t$ may be viewed as the composition of
taking product with Euler class of 
fiber tangent bundle $\tau$ and a pushforward homomorphism. In other words, we
have $c_r(\xi)=p_!(e(\gamma)e(\tau))$, so the Poincare dual of $c_r(\xi)$ is the
image of a fundamental class of nonsingular manifold -- a transverse
intersection of Poincare duals of $e(\gamma)$ and $e(\tau)$.

The classes $c_r^U(\xi)$ and $t(e^U(\gamma))$ also
coincide in complex cobordisms, as shown in \cite{viniti}. But as we show
later, there are (very simple) examples 
when the set $X_1(\xi)$ is nonsingular but its
cobordism class is not equal to $c_1^U(\xi)$. Hence, the construction in
Prop. 2.3 does {\it not} give the resolution of singularities of $X_r(\xi)$ in
general.

We finish this section with just one more 
geometric construction, which works only for $c_1(\xi)$.\\ 

{\bf Proposition 2.4} {\it Let $\det\xi=\Lambda^n \xi$ be a determinant line bundle
for $\xi$. Then $c_1(\xi)=c_1(\det\xi)=e(\det\xi)$.}\\

The easiest way to prove this statement is to use splitting principle: if
$t_1\dots t_n$ are Chern roots of $\xi$, then $c_1(\Lambda^n \xi)=c_1(\xi)=t_1+\dots+t_n$. $\Box$

Since $\det\xi$ is a linear bundle, the Poincare dual to its first Chern class
may be represented by oriented nonsingular submanifold of codimension two. We
obtain another geometric realization of $c_1(\xi)$ which has no analogue for
$c_r(\xi)$ if $r>1$.

\section{Complex cobordisms and the problem of resolution.}

As we said before, the class $c_r(\xi)\in H^{2r}(M^{2m},\m Z)$ may be defined as
Poincare dual to the set $X_r(\xi)$ where the generic section $s$ of bundle
$\op{Hom}(\m C^{n-r+1},\xi)$ is of non-maximal rank. 

In this section we consider the case when $X_r(\xi)$ is a nonsingular manifold for any
generic section $s$. This condition always holds if $n=1$, $r=n$ (this
corresponds to Euler class) or $r=1$ and $m<4$ (\cite{AVGZ}).\\

{\bf Proposition 3.1.} {\it Suppose $X_r(\xi)$ is nonsingular. Then its
normal bundle carries canonical complex structure.}\\

Let $U\subset M$ be a neighbourhood of a point $x\in X_r(\xi)$ such that $\xi$
is trivial over $U$. Then $n-r+1$ sections of $\xi$ over $U$ are determined by a
$C^{\infty}$-map $f:U\to (\m C^n)^{n-r+1}$. Hence $X_r(\xi)$ is the set of
points where vanish $r$ corner minors of the corresponding matrix 
$n\times(n-r+1)$. So $X_r(\xi)$ is locally determined by $r$ complex equations
-- it's an intersection of $r$ nonsingular manifolds of (real) codimension two
with complex structures in their normal bundle. The corresponding complex
structure is preserved under change of trivialization. $\Box$

{\bf Example 3.2.} Let $\xi=\e O(1)\oplus \e O(1)$ be a 2-dimensional
bundle over $M=\m CP^2$, $r=1$. 

Consider two generic holomorphic sections of $\xi$ over $\m CP^2$ -- they have
form $(a_1z_1+a_2z_2+a_3z_3, b_1z_1+b_2z_2+b_3z_3)$ and
$(c_1z_1+c_2z_2+c_3z_3, d_1z_1+d_2z_2+d_3z_3)$, where $z_1,z_2,z_3$ are
homogenous coordinates on $\m CP^2$. Then the condition that sections are
linearly dependent means that
$\frac{a_1z_1+a_2z_2+a_3z_3}{b_1z_1+b_2z_2+b_3z_3}$ = 
$\frac{c_1z_1+c_2z_2+c_3z_3}{d_1z_1+d_2z_2+d_3z_3}$. So $X_1(\xi)$ is a
nonsingular plane curve of degree two.

It is known that Chern classes $c_r^U(\xi)$ 
may also be defined in complex cobordisms --
they're uniquely determined by the same four axioms that define cohomological
Chern classes. It turns out that if $r<n$ then cobordism classes corresponding
to nonsingular $X_r(\xi)$'s are not necessarily equal to $c_r^U(\xi)$.\\

{\bf Proposition 3.3.} {\it Classes $c_r^U(\xi)$ and $[X_r(\xi)]$ are not equal
when $M=\m CP^2$, $r=1$, $\xi=\e O(1)\oplus\e O(1)$.}\\

We'll calculate these classes explicitly. Denote by $t=c_1^U(\e O(1))$ the
cobordism class of the line $\m CP^1\subset \m CP^2$, then, by Whitney sum
formula, $c_1^U(\e O(1)\oplus\e O(1))=c_1^U(\e O(1))+c_1^U(\e O(1))$. 

In our example $X_1(\xi)$ is a plane quadric, so its cobordism class equals to
$e^U(\e O(2))=c_1^U(\e O(2))$. This class is the sum of two classes $c_1^U(\xi)$
in the formal group of geometric cobordisms:
$c_1^U(\nu\oplus\eta)=c_1^U(\nu)+c_1^U(\eta)-[\m
CP^1]c_1^U(\nu)c_1^U(\eta)+\dots$ (\cite{BMN}). So $[X_1(\xi)]=2t-[\m CP^1]t^2$,
which is not $2t$. $\Box$\\

{\bf Problem 3.4.} {\it Construct the classes $A_r(\xi)\in U^{2r}(M)$ satisfying
the following properties:
\begin{enumerate}

{\item Classes $A_r(\xi)$ are functorial.}

{\item $A_r(\xi)=[X_r(\xi)]$ if $X_r(\xi)$ is nonsingular.}

{\item If $X_r(\xi)$ is singular then $A_r(\xi)$ is realized by complex-oriented
map (\cite{KF}) of a manifold resolving singularities of $X_r(\xi)$.}
\end{enumerate}}

Also we expect some more nice properties from $A_r(\xi)$ -- for example, they
could satisfy some Chern classes axioms. 

By functoriality property, class $A_r(\xi)$ is a formal series in
classes $c_1^U(\xi), c_2^U(\xi), \dots$ over ring $U^*(pt)=\m Z[a_1,a_2,\dots],
\deg a_i=-2i$. The augmentation map $U^*(BU)\to H^*(BU,\m Z)$ acts as follows:
all $c_i^U(\xi)$'s map to $c_i(\xi)$'s and $a_i$'s -- to zero. From the second
and third condition it follows that $A_r(\xi)$ has the form $c_r^U(\xi)+\dots$,
where the part "dots" maps to zero under augmentation map.

The geometric constructions from the previous section may easily be
spread to the complex cobordisms. The most efficient construction are classes
$Q_r(\xi)$, which are defined by using the
projectivization of trivial bundle of rank $n-r+1$. In the next sections we 
investigate some of their
properties and show that they satisfy deformed Whitney sum formula.

Hence, the problem of constructing characteristic classes using singularity
cycles makes us vary the system of corresponding axioms. Classes $Q_r(\xi)$
don't only possess deformed Whitney sum formula, they also may be non-trivial in
negative dimension.

\section{The classes $Q_r(\xi)$ and their basic properties.}

In this section we define characteristic classes $Q_r(\xi)\in U^{2r}(M)$ and
establish some of their properties. Recall that $M$ is an oriented manifold and
$\xi$ -- a complex vector bundle over $M$.

Consider the projectivization $M\times \m CP^{n-r}$ of a trivial vector bundle
$\m C^{n-r+1}$ over $M$ and its tautological bundle $\gamma$. Then the section
$s$ determines the section $s_Q$ of a bundle $\op{Hom}(\gamma, p^*\xi)$ (here
$p:M\times \m CP^{n-r}\to M$ is a projection map) by taking the composition
$\gamma\to \m C^{n-r+1}\to p^*\xi$. 

The set $X_r^Q(\xi)$ of points on $M\times\m CP^{n-r}$ where $s_Q$ vanishes is a
nonsingular submanifold, its complex
cobordism class is well-defined and coincides with Euler class
$e^U(\gamma^*\otimes p^*\xi)$.\\

{\bf Definition 4.1.} {\it $Q_r(\xi) = p_! e^U (\gamma^*\otimes p^*\xi)$.}\\

{\bf Theorem 4.2.} {\it \begin{enumerate}

{\item Classes $Q_r(\xi)$ are functorial.}

{\item Classes $Q_r(\xi)$ satisfy the dimension axiom: $Q_r(\xi)=0$ if
$r>n$.}

{\item Classes $Q_r(\xi)$ are normalized:
$Q_n(\xi)=e^U(\xi)=c^U_n(\xi)$.}

{\item If $X_r(\xi)$ is a nonsingular manifold in $M$, then
$Q_r(\xi)=[X_r(\xi)]$.}

{\item The Whitney sum formula fails for $Q_r(\xi)$ so they don't
coincide with $c_r^U(\xi)$.}
\end{enumerate}}

The functoriality property follows from construction. If $r>n$, then $n-r+1$ is
non-positive and class $Q_r(\xi)$ corresponds to an empty map to $M$, so
$Q_r(\xi)$ is zero. 

If $r=n$, then $M\times \m CP^{n-r} = M$, the bundle $\gamma^*$ is trivial and
$p$ is an identity map. So $Q_r(\xi) = p_! e^U(\gamma\otimes p^*\xi) =
e^U(\gamma^*\otimes\xi) = e^U(\xi)$. This proves the third statement.

The class $Q_r(\xi)$ is, by definition, realized by map of a nonsingular
manifold $X_r^Q(\xi)$ to $M$, resolving the singularities of $X_r(\xi)$. If
$X_r(\xi)$ is nonsingular itself, then the corresponding map is identical and
$Q_r(\xi)=[X_r(\xi)]$. 

To prove the last statement it's enough to find a counterexample to Whitney sum
formula. We'll do it in the next section. $\Box$

\section{Examples: bundles over projective spaces.}

{\bf Example 5.1.} Let $\xi$ be an arbitrary linear bundle over $M=\m
CP^1$, $e^U(\xi)=c_1^U(\xi)=u$. We'll calculate $Q_r(\xi)$ for $r=1, 0, -1$.

If $r=1$, then $n-r+1=1$, $M\times\m CP^{n-r} = M = \m CP^1$ and $\gamma$ is
trivial. So $Q_1(\xi)=p_! e^U(\gamma^*\otimes p^*\xi) = e^U(\gamma\otimes\xi) =
e^U(\xi) = u$. 

If $r=0$, then $M\times \m CP^{n-r}$ = $\m CP^1\times \m CP^1$, $p$ is a
projection map to first $\m CP^1$. Denote $e^U(\gamma^*)$ by $t$. Then the class
$e^U(\gamma^*\otimes p^*\xi)$ has the form $t+u-[\m CP^1]tu+\dots$ (\cite{BMN}),
where every monomial in ($\dots$) is divisible by either $t^2$ or $u^2$.

The pushforward homomorphism $p_!$, corresponding to a map $M\times \m CP^k\to
M$ is a $U^*(M)$-module homomorphism given by formula 
$p_!(t^{i})=[\m CP^{k-i}]$ if $i\leqslant k$ and $p_!(t^i)=0$ if $i>k$. In our
case $p_!(t+u-[\m CP^1]tu+\dots) = 1 - [\m CP^1]u + [\m CP^1]u +\dots$ because
$u^2=0$ in $U^*(\m CP^1)$. So $Q_0(\xi)=1$ 
for every linear bundle $\xi$ over $\m
CP^1$. 

Finally, if $r=-1$ then $n-r+1=3$ and $p:\m CP^1\times\m CP^2\to \m CP^1$ is
again a projection map. In this case $p_!(t^2)=1$, $p_!(t)=[\m CP^1]$,
$p_!(1)=[\m CP^2]$ and $Q_{-1}(\xi)=p_! e^U(\gamma^*\otimes p^*\xi) = p_! (t+u-[\m
CP^1]tu+\dots) = [\m CP^1]+([\m CP^2]-[\m CP^1]^2)u \ne 0\in U^{-2}(\m CP^1)$.
We see that classes $Q_r(\xi)$ and $c_r^U(\xi)$ are not equal because
$c_r^U(\xi)=0$ if $r<0$.  

{\bf Example 5.2.} Let $\xi$ be a bundle of rank $k>1$ over $\m CP^1$. We'll
show that again $Q_0(\xi)$ is $1$.

Let $t_1\dots t_k$ be Chern roots of $\xi$, then $Q_0(\xi)$ = $p_!
e^U(\gamma^*\otimes p^*\xi) = p_! (F(t,t_1)\dots F(t,t_k)) = p_!((t+t_1-[\m
CP^1]tt_1+\dots)\dots(t+t_k-[\m CP^1]tt_k+\dots))$, where $p:\m CP^1\times \m
CP^k\to \m CP^1$ is a projection map. 

In our case $p_!(t^k)=1$, $p_!(t^{k-1})=[\m CP^1]$ and $t_i\cdot t_j=0$ in
$U^*(\m CP^1)$. So $Q_0(\xi)=p_!(t^k + t^{k-1} (t_1+\ldots+t_k) - [\m
CP^1]t^k(t_1+\ldots+t_k)) = 1 - [\m CP^1](t_1+\ldots+t_k) + [\m
CP^1](t_1+\ldots+t_k) +\dots = 1$.

We see that in case of the bundle over $\m CP^1$ classes $Q_0(\xi)$ and
$Q_1(\xi)$ coincide with $c_0^U(\xi)$ and $c_1^U(\xi)$ respectively. 

{\bf Example 5.3.} Let us show that $Q_0(\e O(1))=1+u^2([\m CP^1]^2-[\m CP^2])$
in $U^0(\m CP^2)$. If $p:\m CP^1\times\m CP^2\to\m CP^2$ 
is a projection, $t=e^U(\gamma^*)$, $u=e^U(p^* \e O(1))$, then $Q_0(\e O(1))$ = $p_!
e^U(\gamma^*\otimes p^*\xi)$ = $p_! (t+u-[\m CP^1]tu+([\m CP^1]^2-[\m
CP^2])tu^2+\dots)$, where any monomial in $(\dots)$ divides either by $t^2$ or
$u^3$. So $Q_0(\xi)$ = $1+[\m CP^1]u - [\m CP^1]u+([\m CP^1]^2-[\m CP^2])u^2$ =
$1 + ([\m CP^1]^2-[\m CP^2])u^2$ in $U^0(\m CP^2)$.

We may now give the example of bundles $\xi$ and $\eta$ such that 
classes $Q_r(\xi\oplus\eta)$, $Q_i(\xi)$ and $Q_j(\eta)$ don't satisfy
Whitney sum formula. Consider $M=\m CP^2$, $\xi = \eta = \e O(1)$. Then
$Q_1(\e O(1)\oplus \e O(1))$ is equal to cobordism class of a nonsingular
quadric in $\m CP^2$, which is $2u-[\m CP^1]u^2$, $u=e^U(\e O(1))$. 
If we assume that 
Whitney sum formula holds for classes $Q_r(\cdot)$, 
then $Q_1(\e O(1)\oplus\e O(1)) = 2 Q_0 (\e O(1)) Q_1(\e O(1)) = 2u\ne 2u-[\m
CP^1] u^2$. This finishes the proof of Th. 4.2.

As we said before, in general 
the classes $Q_r(\xi)$ are formal series in $c_k^U(\xi)$ with
coefficients in $U^*(pt)$ (this follows from functoriality property). 
For example, if $\xi$
is a 2-dimensional bundle, then the class $Q_1(\xi)$ has the form
$$
Q_1(\xi) = c_1^U(\xi)+[\m CP^1]c_2^U(\xi) + ([\m CP^1]^2-[\m
CP^2])c_1^U(\xi)c_2^U(\xi)+\dots
$$
(which may be computed by splitting principle).

\section{An explicit sum formula for $Q_r(\xi)$.}

In this section we calculate deformed sum formula for classes $Q_r(\xi)$.

Let $\xi$ be a complex vector bundle of rank $k$ and $t_1\dots t_k$ -- its Chern
roots. Consider the product $F(t,t_1)\cdot\ldots\cdot F(t,t_k)$ where $F(u,v)$
is a formal group in geometric cobordisms. Then we can write
$$
\prod\limits_{i=1}^{k} F(t,t_i) = \sum\limits_{j=0}^{\infty} t^j
\Phi_{k-j}(t_1,\dots, t_k),
$$
where $\Phi_r(t_1,\dots,t_k)$ is a symmetric formal series in $t_1,\dots,t_k$
with coefficients in $U^*(pt)$. It is easy to see that
$\Phi_r(t_1,\dots,t_k)=\Phi_r(\xi)$ is a homogenous characteristic class of
dimension $2r$. 

Classes $\Phi_r(\xi)$ naturally arise in the following well-known construction:
if $\xi$ is a complex vector bundle of rank $n$ 
over base $X$, then one can consider 
the bundle $\xi\otimes\gamma$ over $X\times\m CP^{\infty}$ ($\gamma$ is a
tautological linear bundle over $\m CP^{\infty}$). The ring of complex
cobordisms $U^*(X\times\m CP^{\infty})$ is isomorphic to ring of formal series
$U^*(X)[[u]]$, where $\deg u=2$. This means that 
the class $c_n(\xi\otimes\gamma)$ may be written in the form
$\Phi_n(\xi)+u\Phi_{n-1}(\xi)+u^2\Phi_{n-2}(\xi)+\dots$.\\

{\bf Proposition 6.1.} {\it Classes $\Phi_r(\xi)$ satisfy the following
properties:
\begin{enumerate}
{\item $\Phi_r(\xi)=0$ if $r>n$.}
{\item Classes $\Phi_r(\xi)$ are functorial.}
{\item Classes $\Phi_r(\xi)$ satisfy Whitney sum formula
$\Phi_r(\xi\oplus\eta)$ = 
$\sum\limits_{j=-\infty}^{\infty}\Phi_j(\xi)\cdot\Phi_{r-j}(\eta)$ (which is
well-defined by property 1).}
{\item Class $\Phi_n(\xi)$ is equal to Euler class (top Chern class) $e^U(\xi)$ =
$c_n^U(\xi)$ of bundle $\xi$.}
\end{enumerate}}

These properties are immediate from the definitions; for example, the Whitney
sum formula follows from identity $(\xi\oplus\eta)\otimes\gamma =
(\xi\otimes\gamma)\oplus(\eta\otimes\gamma)$. 

By using upper-triangular change of variables one can express $Q_r(\xi)$ through
$\Phi_r(\xi)$ (and vice versa).\\

{\bf Proposition 6.2.} 
{\it $Q_r(\xi)=\sum\limits_{k=0}^{n-r}\Phi_{r+k}(\xi)[\m CP^k]$.}\\

As we said before, the pushforward homomorphism 
$p_!:H^*(M\times\m CP^{k})\to H^*(M)$ 
is a $U^*(pt)$-module homomorphism given by formula $p_!(t^i)=[\m CP^{k-i}]$ if
$i\leqslant k$, $p_!(t^i)=0$ otherwise. So one can obtain $Q_r(\xi)$ by taking
$c_n(\gamma^*\otimes p^*\xi)$  = $\Phi_n(\xi)+t\Phi_{n-1}(\xi)+\dots$ and substitution $t^i\to [\m CP^{n-r-i}]$. $\Box$

The transition matrix $\Phi_r(\xi)\to Q_r(\xi)$ is upper-triangular; its $i$-th
diagonal is containing only $[\m CP^i]$'s. The inverse matrix is also
upper-triangular: its $i$-th diagonal consists of elements $[M^i]\in
U^{-2i}(pt)$, where generating function $1+\sum[M^i]x^i$ is equal to
$\frac1{1+\sum[\m CP^i]x^i}$.\\

{\bf Theorem 6.3.} {\it Consider classes $[M^i]\in
U^{-2i}(pt)$ with generating function 
$1+\sum[M^i]x^i$ = $\frac1{1+\sum[\m CP^i]x^i}$. Then for any complex vector
bundles $\xi$, $\eta$ and any $r$ we have
$$
Q_r(\xi\oplus\eta)=\sum\limits_{i,j,k,l\in\m Z}
Q_{l+i}(\xi)Q_{r-l+k+j}[M^i][M^j][\m CP^k].
$$
(if $i<0$, then $[M^i]=[\m CP^i]=0$).}

The proof is straightforward: $Q_r(\xi\oplus\eta)$ = $\sum\limits_{k\in\m Z}
\Phi_{r+k}(\xi\oplus\eta)[\m CP^k]$ = $\sum\limits_{k,l\in\m
Z}\Phi_l(\xi)\Phi_{r-l+k}(\eta)[\m CP^k] = \sum\limits_{k,l,i,j\in \m
Z}Q_{l+i}(\xi)Q_{r-l+k+j}[M^j][\m CP^k]$. $\Box$

If we set $[M^i]=[\m CP^i]=0$ for all positive $i$ (this corresponds to
augmentation map $U^*(pt)\to H^*(pt)$), then our sum formula turns into 
standart Whitney sum formula for classes $c_r(\xi)$ -- as it should be.\\

{\bf Theorem 6.4.} {\it Classes $Q_r(\xi)\in U^{2r}(X)$ of a complex vector
bundle bundle $\xi$ of rank $n$ over $X$ are uniquely determined by four axioms:
\begin{enumerate}
{\item Classes $Q_r(\xi)$ are functorial.}
{\item $Q_r(\xi)=0$ if $r>n$.}
{\item Classes $Q_r(\xi)$ satisfy the sum formula 
$$
Q_r(\xi\oplus\eta)=\sum\limits_{i,j,k,l\in\m
Z}Q_{l+i}(\xi)Q_{r-l+k+j}(\eta)[M^i][M^j][\m CP^k]
$$}
{\item $Q_n(\xi)$ = $e^U(\xi)$ = $c_n^U(\xi)$.}
\end{enumerate}}

The proof is analogous to proof of uniqueness of cohomological Chern classes.
Axiom 4 may be reformulated in a more classical way: the class $Q_1(\e O(1))$,
where $\e O(1)$ is a canonical bundle over $\m CP^1$, 
is equal to fundamental class $u\in U^2(\m CP^1)$.

\section{Classes $P_r(\xi)$ and $IH$-small resolutions.}

The definition of classes $P_r(\xi)$ follows the construction from section 1.
Given a generic section $s:\m C^{n-r+1}\to \xi$ (we denote by $\m C$ the trivial
linear bundle) one can consider the dual section $s^*:\xi^*\to\m C^{n-r+1}$. If
$\dim\ker s=l$ in some point $x\in M$, then $\dim\ker s^* = l+r-1$. 

Let $\gamma$ be an $r$-dimensional tautological bundle over $Gr_r(\xi^*)$ --
$r$-th grassmanization of $\xi^*$, $p:Gr_r(\xi^*)\to M$ -- a projection map.
Then section $s^*$ determines the section $\gamma\to p^*\xi^*\to \m C^{n-r+1}$
of a bundle $\op{Hom}(\gamma,\m C^{n-r+1})$ over $Gr_r(\xi^*)$. Let $X_r^P(\xi)$
be a zero set of this new section. If $\dim\ker s=l$ in $x\in M$,  then
$p^{-1}(x)\cap X_r^P(\xi) = Gr_r(\m C^{l+r-1})$. The set $X_r^P(\xi)$ is a
nonsingular manifold and its cobordism class is equal to $e^U(\op{Hom}(\gamma,\m
C^{n-r+1}))$ = $e^U((\gamma^*)^{n-r+1})$.\\

{\bf Definition 7.1.} $P_r(\xi) = p_! e^U((\gamma^*)^{n-r+1})$.\\

{\bf Theorem 7.2.} {\it
\begin{enumerate}
{\item Classes $P_r(\xi)$ are functorial.}
{\item $P_r(\xi)=0$ if $r>n$.}
{\item The classes $P_r(\xi)$ satisfy the normalization axiom:
$P_n(\xi)=e^U(\xi)$.}
{\item If set $X_r(\xi)$ is nonsingular then $P_r(\xi)=[X_r(\xi)]$.}
{\item The classes $P_r(\xi)$ and $c_r^U(\xi)$ are not equal.}
\end{enumerate}}

The proofs are straightforward from definitions (except the last statement which
follows from theorems below). 

We see that there exist at least two resolutions of singular manifold
$X_r(\xi)$, corresponding to $Q_r(\xi)$ and $P_r(\xi)$ respectively. The classes
$P_r(\xi)$ and $Q_r(\xi)$ differ -- but only in high dimensions.\\

{\bf Theorem 7.3.} {\it
\begin{enumerate}
{\item If $\dim\xi=2$, the class $P_1(\xi)-Q_1(\xi)$ vanishes if $\dim_{\m R}
M\leqslant 8$ and may be nonzero if $\dim_{\m R} M=10$ (recall that $M$ is
even-dimensional).}
{\item If $\dim\xi=2$, the class $P_2(\xi)-Q_2(\xi)$ vanishes if $\dim_{\m R}
M\leqslant 10$ and may be nonzero if $\dim_{\m R} M=12$.}
\end{enumerate}}

This is proved by direct computation using generalized Riemann-Roch theorem for
Chern-Dold character (\cite{buch}). Chern-Dold character is a
multiplicative transform $\op{ch}_U:U^*(X)\to H^*(X,U^*(pt))$, which is an
isomorphism for $X=BU(n)$. So if one wants to compare $P_r(\xi)$ and $Q_r(\xi)$,
then it is better 
to compare their transformations 
$\op{ch}_U P_r(\xi)$ and $\op{ch}_U Q_r(\xi)$ -- the
result should be the same, but the computations will be much more easy. If
$p:X\to Y$ is a smooth bundle, $\tau$ -- its fiber tangent bundle carrying a
complex structure, then (according to generalized Riemann-Roch theorem)
$$
\op{ch}_U p_!^U(x) = p_!^H (\op{ch}_U(x)\cdot T(\tau)),
$$
where $p_!^U$ and $p_!^H$ are pushrorward homomorphisms in cobordisms and
cohomology and $T(\tau)\in H^*(X,U^*(pt))$ is generalized Todd class. $T(\tau)$
is multiplicative: $T(\eta\oplus\xi)=T(\eta)T(\xi)$ and if $\gamma$ is linear,
then $T(\gamma)=\frac{c_1(\gamma)}{g^{-1}(c_1(\gamma))}$ (see \cite{buch}),
where $g^{-1}$ is an exponential of formal group of complex cobordisms.
The homomorphism $p_!^H$ is easier than $p_!^U$, so we can compare
$\op{ch}_U Q_r(\xi)$ and $\op{ch}_U P_r(\xi)$ explicitly.

Let $\dim\xi=2$ and $r=1$. If $\dim M\leqslant 6$ then necessarily $Q_1(\xi)$ = $P_1(\xi)$ because
$X_1(\xi)$ is nonsingular. 
The beginning of the row corresponding to $Q_1(\xi)$ and
$P_1(\xi)$ has the form
$$
c_1^U(\xi) - [\m CP^1]c_2^U(\xi) + ([\m CP^1]^2-[\m
CP^2])c_1^U(\xi)c_2^U(\xi)+\dots
$$

If $\dim M=8$ or $10$, then $X_1(\xi)$ may be singular, but still there is no
difference between $Q_1(\xi)$ and $P_1(\xi)$. This is not surprising: the
manifolds $X_1^P(\xi)$ and $X_1^Q(\xi)$ are "made of equal parts": the preimages
of point $x\in X_1(\xi)$ where $\dim\ker s=l$ are diffeomorphic 
to $\m CP^{l-1}$ in $X_1^P(\xi)$ and $X_1^Q(\xi)$. 

Not let $M$ be a complex variety and $\xi$ a holomorphic vector bundle,
$\dim_{\m C} M=4$, $\dim\xi=2$. Then $X_1(\xi)$ is nonsingular in points of 
the complement
to several isolated points, and in these points 
$X_1(\xi)$ is locally given by equation $xy-zt=0$. Then $X_1^Q(\xi)$ and
$X_1^P(\xi)$ correspond to two different resolutions of isolated quadratic
singularity.

This corresponds to the following situation: the resolution $f:X\to Y$ of a
singular variety $Y$ is called $IH$-small if $\op{codim}_{\m C} Y_i>2i$, where 
$Y_i=y\in Y\vert \dim_{\m C} f^{-1}(y) = i$. Most singular 
varieties don't possess $IH$-small resolutions. As proved by B.Totaro
(\cite{totaro}), two different $IH$-small resolutions of one singular variety
have equal elliptic genera. So even if $P_1(\xi)\ne Q_1(\xi)$,
the classes $[X_1^P(\xi)]$ and $[X_1^Q(\xi)]$ have the same elliptic genus.

Now let $\dim\xi=3$ and $r=2$. If $\dim M\leqslant 10$ then $X_2(\xi)$ is again
nonsingular and $Q_2(\xi)=P_2(\xi)$. The beginning of the row corresponding to
$Q_2(\xi)$ and $P_2(\xi)$ has the form
$$
c_2^U(\xi)-2[\m CP^1]c_3^U(\xi) + ([\m CP^1]^2-[\m
CP^2])c_1^U(\xi)c_3^U(\xi)+\dots
$$

The preimage in $X_2^Q(\xi)$ of point $x\in X_2(\xi)$ where $\dim\ker s=l$ is
$\m CP^1$ and in $X_2^P(\xi)$ it is $Gr_2(\m C^3) = \m CP^2$. So Euler
characteristics of $X_2^P(\xi)$ and $X_2^Q(\xi)$ are not equal if $X_2(\xi)$ is
singular. The corresponding cobordism classes $Q_2(\xi)$ and $P_2(\xi)$ also
differ.\\

{\bf Corollary 7.4.} {\it Two resolutions of singular variety 
$X_r(\xi)$ corresponding to $Q_r(\xi)$ and $P_r(\xi)$ are not equivalent.}\\

If $r>1$, this is obvious because Euler characteristics of corresponding
manifolds differ; the case $r=1$
follows from the fact that classes $P_1(\xi)$ and $Q_1(\xi)$ aren't
equal. $\Box$

\section{The class $D_1(\xi)$.} 

Finally, in this section 
we consider the characteristic class $D_1(\xi) = c_1^U(\det\xi)\in U^2(M)$ of 
$\det\xi=\Lambda^n\xi$. It is functorial and equal to $[X_1(\xi)]$ if $X_1(\xi)$
is nonsingular. But our last condition (resolution of singularities of
$X_1(\xi)$ if it's singular) fails for $D_1(\xi)$ (example: $\e O(1)\oplus\e
O(-1)$ over $\m CP^2$). But if $X_1(\xi)$ is nonsingular, the class $D_1(\xi)$
is much easier to compute than $Q_1(\xi)$ and (of course) $P_1(\xi)$.\\

{\bf Proposition 8.1.} {\it We have 
$$
\op{ch}_U D_1(\xi) = \sum\limits_{i=0}^{\infty} [N^{2i}]c_1(\xi)^{i+1},
$$
where elements $[N^{2i}]\in U^*(pt)\otimes\m Q$ satisfy $s_{(i)}([N^{2i}])=1$,
$s_{\omega}([N^{2i}])=0$ if $(\omega)\ne (i)$, $s_{\omega}$ are
Landweber-Novikov operations in complex cobordism.}\\

We use splitting principle for the one last time. If
$\xi=\eta_1\oplus\dots\oplus\eta_n$ is a direct sum of linear bundles,
$c_1^U(\eta_i)=t_i$, then $c_!^U(\det\xi) = g^{-1}(g(t_1)+\ldots+g(t_k))$, where
$g(x)=\sum\limits_{i=0}^{\infty} \frac{[\m CP^i]}{i+1} x^{i+1}$ is a logarithm
of formal group of complex cobordisms. 

The Chern-Dold character is a ring homomorphism, 
so $\op{ch}_U g^{-1} (g(t_1)+\ldots+g(t_n)) = g^{-1}( g(\op{ch}_U t_1)
+\ldots+g(\op{ch}_U t_k))$. But, as we said before, $\op{ch}_U(x)$ regarded as
formal series, coincides with $g^{-1}(x)$ (\cite{buch}). So if $\op{ch}_U t_i =
g^{-1}(s_i)$, then $\op{ch}_U g^{-1}(g(t_1)+\ldots+g(t_n)) =
g^{-1}(s_1+\ldots+s_n) = g^{-1}(c_1(\xi))$. $\Box$

\end{document}